\documentclass[journal]{IEEEtran}
\usepackage[short]{optidef}
\usepackage{graphicx, float}
\usepackage{amsmath}
\usepackage{siunitx}
\usepackage{amssymb}
\usepackage{amsthm}
\usepackage{booktabs}
\usepackage{amsfonts}
\usepackage{arydshln}
\usepackage{mathtools}
\usepackage{cite}
\usepackage[subtle]{savetrees}
\usepackage{soul}
\usepackage{color} % load the package
\usepackage[linesnumbered, ruled]{algorithm2e}
\newtheorem{theorem}{Theorem}[section]

\newtheorem{lemma}[theorem]{Lemma}
\newtheorem*{remark}{Remark}

\ifCLASSOPTIONcompsoc
    \usepackage[caption=true, font=normalsize, labelfont=sf, textfont=sf]{subfig}
\else
\usepackage[caption=false, font=footnotesize]{subfig}
\fi
\IEEEoverridecommandlockouts  

\title{%Need New Title
 Guaranteeing a physically realizable battery dispatch without charge-discharge complementarity constraints
}
\author{Nawaf~Nazir$^\ast$ and Mads~Almassalkhi$^\dagger$
\thanks{$^{\ast}$N. Nazir is affiliated with Pacific Northwest National Lab, Richland, WA 99354.      $^{\dagger}$M. Almassalkhi is affiliated with the Department of Electrical and Biomedical Engineering, The University of Vermont, Burlington, VT 05405, USA.

Support from the U.S. Department of Energy’s Advanced Research Projects Agency—Energy (ARPA-E) Award DE-AR0000694 is gratefully acknowledged. M. Almassalkhi was supported in part by the National Science Foundation (NSF) Award ECCS-2047306.}}%
\date{}
\begin{document}
\maketitle
\begin{abstract}\label{sec:abstract}
The non-convex complementarity constraints present a fundamental computational challenge in energy constrained optimization problems. In this work, we present a new, linear, and robust battery optimization formulation that sidesteps the need for battery complementarity constraints and integers and prove analytically that the formulation guarantees that all energy constraints are satisfied which ensures that the optimized battery dispatch is physically realizable.
%This has the potential to greatly reduce the  of energy constrained optimization problems. 
%Through analysis, we show that this new formulation is  under practical conditions and
In addition, we bound the worst-case model mismatch and discuss conservativeness. Simulation results further illustrate the effectiveness of this approach.
% In this work, we develop a formulation that avoids simultaneous charging and discharging of batteries without the use of binary variables. 
%Mathematical analysis shows that this formulation holds under practical conditions. 
%Simulation results further illustrate the effectiveness of the approach.
\end{abstract}
\begin{IEEEkeywords}
Energy storage, battery, simultaneous charging and discharging, complementarity constraint, model predictive control.
\end{IEEEkeywords}
\section{Introduction}\label{sec:introduction}
Due to the increasing penetration of variable renewable generation, widespread concerns over reliability of power systems are being raised. Deploying battery storage systems is  widely considered a solution to improve grid operations and reliability~\cite{almassalkhi2014model}. However, optimizing battery storage requires being cognizant of the dynamics of the state of charge (SoC), limits on SoC, limits on the rate of change of the SoC (i.e., power input/output), and the physical operating modes of the battery: it can either charge (i.e., consume energy) or discharge (i.e., produce energy), but not both at the same time. Previous works in literature have employed binary variables (i.e., $1=$charging and $0=$discharging) in the optimization to overcome this issue~\cite{hu2008global}. However, solving general mixed-integer programs (MIPs) is computationally challenging.

To sidestep the MIP challenge, several works in literature have proposed different battery models as a way to overcome the non-convex complementarity constraints~\cite{li2015sufficient,aaslid2020non,garifi2020convex,nazir2020optimal,almassalkhi2014model}. {\color{black} In~\cite{aaslid2020non}, the battery physics are modeled in the nonlinear current-voltage variable space, which gives a non-convex, but continuous formulation.} The authors in~\cite{li2015sufficient} relaxed the complementarity constraints in a bulk transmission economic dispatch problem. Then, through KKT analysis, they provide sufficient conditions under which the relaxation holds. However, these conditions do apply under negative locational marginal prices (LMPs). Similar to the work in~\cite{li2015sufficient}, the authors in~\cite{garifi2020convex,nazir2020optimal} extend the formulation to distribution networks and provide methods to avoid simultaneous charging and discharging by modifying the objective function. However, these methods do not hold under high renewable penetration, specifically under reverse power flow. 
%The authors in~\cite{hari2018hierarchical} introduce auxiliary variables, which are then related to the original charging and discharging variables through convex relaxations. They empirically show that the relaxation is tight, but do not provide any theoretical guarantees. 
The work in~\cite{almassalkhi2014model} quantifies the effects of simultaneous charging and discharging and also provides a heuristic approach to avoid this phenomenon. Recently~\cite{arroyo2020use} reiterated that many of these approaches fail in practical settings and engender simultaneous charging and discharging.

From the above literature, battery models based on relaxing complementarity constraints fail under practical conditions which then leads to violation of battery SoC constraints when the said models are employed in optimization problems. This is particularly true with reverse power flow from increasing number of vehicle to grid (V2G) systems. Hence, there is a need for models that respect the SoC constraints without having to resort to non-linear complementarity constraints.  Another factor is a shift towards real-time control of power systems~\cite{tang2017real}, which motivates a need to avoid mixed-integer formulations. To tackle this critical problem, in this work, we propose a method that respects the battery SoC constraints under practical conditions and, at the same time, avoids the need for non-linear complementarity constraints or binary variables. We augment the battery model with a linear term that utilizes a simplified battery model using only the net battery power exchanges. This simplified linear term results in tightening of the SoC upper limit in the battery model.
%To account for the modeling errors introduced by this term (due to the non-unity charging and discharging efficiencies), we tighten the SoC lower limit in the battery model.
The contribution is a new linear energy storage dispatch formulation whose optimal solution predicts a physically realizable dispatch, i.e., a sequence of power set-points whose resulting SoC trajectory respects the actual SoC limits. \textcolor{black}{ That is, the authors' definition of \textit{physically realizable} refers to the state (or system's output) being achievable rather than an implementation (i.e., systems input sequence) being achievable}. We provide analysis that proves the feasibility of this technique and also provide bounds on its conservativeness with regards to the tightening of the SoC limits. Bounds on the worst-case tightening of SoC limits can be calculated \textit{a priori} based on parameters such as the optimization horizon and time-step and the battery specs.

In the rest of the manuscript we develop a novel, linear formulation of the optimal battery dispatch problem that respects the SoC limits without using complementarity constraints. We analyze the approach and provide simulation results that demonstrate its effectiveness.
% The rest of the manuscript is organised as follows: Section~\ref{sec:stan_model} presents the standard battery SoC model. In section~\ref{sec:simp_model} we introduce the simplified battery model that avoids simultaneous charging and discharging and quantify the error introduced due to the simplified model, whereas in section~\ref{sec:prob_form} we provide the optimal battery dispatch formulation. Section~\ref{sec:simulations} provides simulation results that illustrate the effectiveness of the approach and finally, the manuscript is concluded in Section~\ref{sec:conclusion}, which also outlines directions for future work.
\section{Standard battery model}\label{sec:stan_model}
Consider a battery with SoC at (discrete) time-step $k$, $E[k]\in [0,E_{\text{max}}]$, where each time-step represents a duration $\Delta t>0$. The battery also has charging and discharging inputs that can be applied over time-step $k$ defined as $P_{\text{c}}[k], P_{\text{d}}[k]\in [0,P_{\text{max}}]$, respectively, and charging and discharging efficiencies $\eta_{\text{c}},\eta_{\text{d}}\in (0,1]$, respectively. In addition, the battery can either charge or discharge but not both at time $k$, which yields non-convex complementarity condition $P_{\text{c}}[k] P_{\text{d}}[k]= 0$.  Then, starting with a given initial SoC $E_0$ and a sequence of inputs over period $\mathcal{T}=\{0,1,\hdots,T-1\}$, the battery SoC dynamics evolve along a admissible trajectory described by the following set of equalities and inequalities:
\begin{subequations}\label{eq:stan_model}
\begin{align}
    E[k+1]=&E[k]+\Delta t\eta_{\text{c}}P_{\text{c}}[k]-\frac{\Delta t}{\eta_{\text{d}}}P_{\text{d}}[k], \qquad \forall k\in \mathcal{T}\label{eq:SM_1}\\
    E[0]=&E_0\label{eq:SM_2}\\
    0\le P_{\text{c}}[k]\le& P_{\text{max}}, \qquad \forall k\in \mathcal{T}\label{eq:SM_3}\\
    0\le P_{\text{d}}[k]\le& P_{\text{max}}, \qquad \forall k \in \mathcal{T}\label{eq:SM_4}\\
    0\le E[k+1]\le& E_{\text{max}}, \qquad \forall k\in \mathcal{T}\label{eq:SM_5}\\
     P_{\text{c}}[k] P_{\text{d}}[k]=& 0 \qquad \forall k\in \mathcal{T}. \label{eq:SM_cc}
    %z[k]\in &\{0,1\}, \qquad \forall k \in \mathcal{T}\label{eq:SM_6}
\end{align}
\end{subequations}
%Note that the complementarity condition on charging and discharging is captured with the binary variable $z[k]$ in the right-hand sides of~\eqref{eq:SM_3} and~\eqref{eq:SM_4} and implies that if $z[k]=1 \implies P_{\text{d}}[k]=0$ and if $z[k]=0\implies P_{\text{c}}[k]=0$.
The resulting SoC trajectory 
%over time-steps $k \in \mathcal{T}$ in~\eqref{eq:SM_1} 
can be expressed as
\begin{align}\label{eq:P1}
    \mathbf{E(P_{\text{c}},P_{\text{d}})}=\mathbf{1}_TE_0+\eta_{\text{c}}\mathbf{AP_{\text{c}}}-\frac{1}{\eta_{\text{d}}}\mathbf{AP_{\text{d}}},
\end{align}
where $\mathbf{E}=\text{col}\{E[k+1]\}_{k\in \mathcal{T}}, \mathbf{P_{\text{c}}}=\text{col}\{P_{\text{c}}[k]\}_{k\in \mathcal{T}}$, and $\mathbf{P_{\text{d}}}=\text{col}\{P_{\text{d}}[k]\}_{k\in \mathcal{T}}$ and %$\mathbf{P_{\text{c}}.P_{\text{d}}=0}$, 
$\mathbf{A}$ is a lower triangular matrix that relates the input at time $k$ to $\Delta t E[l]$ at time $l\ge k$, and $\mathbf{1}_T:=[1,\hdots,1]^\top\in \mathbb{R}^{T}$.

Prior work has relaxed the battery model in~\eqref{eq:stan_model} by removing  complementarity condition~\eqref{eq:SM_cc}~\cite{nazir2020optimal, almassalkhi2014model,li2015sufficient,garifi2020convex}. Thus, relaxed models allow \textit{simultaneous charging and discharging}, i.e., $P_{\text{c}}[k] P_{\text{d}}[k] \ge 0$.

In the next section, we present the relaxed model and a simplified 1-input model that considers only the net-charging input, i.e., $\mathbf{P}_{\text{b}}:=\mathbf{P}_{\text{c}}-\mathbf{P}_{\text{d}}$. Then, we analytically show how these two models together present necessary and sufficient bounds on the SoC in the standard model. This informs a novel battery optimization formulation that is both convex and whose optimal open-loop dispatch schedule is guaranteed to be physically realizable.
%overcomes the shortcomings of the standard model in~\eqref{eq:stan_model}.

\section{Relaxed and simplified battery models}\label{sec:simp_model}
\subsection{Relaxed model}
The relaxed model is obtained by removing~\eqref{eq:SM_cc} from the standard battery model~\eqref{eq:stan_model} and is, thus, convex. %to obtain a linear model. 
It defines a relaxed SoC trajectory $\mathbf{E}^\text{r}:=\text{col}\{E^\text{r}[k]\}_{k=1}^T$ that $\forall k\in \mathcal{T}$ satisfies
%set of equations and constraints:

\begin{subequations}\label{eq:relax_model}
\begin{align}
    E^{r}[k+1]=&E^{r}[k]+\Delta t\eta_{\text{c}}P^{r}_{\text{c}}[k]-\frac{\Delta t}{\eta_{\text{d}}}P^{r}_{\text{d}}[k], \qquad \label{eq:RM_1}\\
    E^{r}[0]=&E_0\label{eq:RM_2}\\
    0\le P^{r}_{\text{c}}[k]\le& P_{\text{max}}, \qquad \label{eq:RM_3}\\
    0\le P^{r}_{\text{d}}[k]\le& P_{\text{max}}, \qquad \label{eq:RM_4}\\
    0\le E^{r}[k+1]\le& E_{\text{max}}, \qquad \label{eq:RM_6}
    %&P^{r}_{\text{c}}[k]+P^{r}_{\text{d}}[k] \le P_{\text{max}}, \qquad \forall k\in \mathcal{T}\label{eq:RM_5}.
\end{align}
\end{subequations}
Note that with complementarity conditions relaxed in~\eqref{eq:relax_model} we have new inputs $P^{r}_{\text{c}}[k]$ and $P^{r}_{\text{d}}[k]$ that are different from those in~\eqref{eq:stan_model}. The variables are related
\begin{align}\label{eq:P_subs}
    \mathbf{P}_{\text{c}}=\max\{\mathbf{0,P}^{r}_{\text{c}}-\mathbf{P}^{r}_{\text{d}}\},\quad \mathbf{P}_{\text{d}}=\max\{\mathbf{0,-(P}^{r}_{\text{c}}-\mathbf{P}^{r}_{\text{d}})\},
\end{align}
which implies that $\mathbf{P}_{\text{c}} - \mathbf{P}_{\text{d}}= \mathbf{P}_{\text{c}}^r-\mathbf{P}_{\text{d}}^r$. 
The relaxed model's SoC trajectory is then defined by
\begin{align}\label{eq:P2}
    \mathbf{E}^{r}(\mathbf{P}^{r}_{\text{c}},\mathbf{P}^{r}_{\text{d}})=\mathbf{1}_T E_0+\eta_{\text{c}}\mathbf{AP}^{r}_{\text{c}}-\frac{1}{\eta_{\text{d}}}\mathbf{AP}^{r}_{\text{d}}.
\end{align}

\subsection{Simplified 1-input model}
For the simplified battery model, we approximate the battery efficiencies, $\eta_{\text{c}}$ and $\frac{1}{\eta_{\text{d}}}$ by a single net-charge efficiency $\eta \in [\eta_{\text{c}}, \frac{1}{\eta_{\text{d}}}]$ and replace the two inputs in~\eqref{eq:stan_model} $P_{\text{c}}[k]$ and $P_{\text{d}}[k]$ with a single net-charging input $P_{\text{b}}[k]=P_{\text{c}}[k]-P_{\text{d}}[k]\in [-P_{\text{max}},P_{\text{max}}]$, which yields the simplified 1-input model:
\begin{subequations}
\begin{align}
    {E}^s[k+1]:=&{E}^s[k]+\eta \Delta t P_{\text{b}}[k], \qquad \forall k\in \mathcal{T}\\
    {E}^s[0]=&E_0\\
    -P_{\text{max}}\le  P_{\text{b}}[k]\le& P_{\text{max}} \qquad \forall k\in \mathcal{T}\\
    0\le  {E}^s[k+1]\le& E_{\text{max}}. \qquad \forall k\in \mathcal{T}
\end{align}
\end{subequations}
 The simplified model's SoC trajectory is then
\begin{align}\label{eq:P3}
    \mathbf{E}^{\text{s}}(\mathbf{P}_{\text{b}})=\mathbf{1}_T E_0 + \eta \mathbf{A}\mathbf{P}_{\text{b}}.
\end{align}

\subsection{Analyzing model mismatch}\label{sec:3_c}
Clearly, by relaxing complementarity conditions and simplifying efficiencies, the corresponding open-loop SoC trajectories may not match the actual trajectory in~\eqref{eq:stan_model}. However, we will next show that the models are ordered in that $\mathbf{E}^{r}\le \mathbf{E}\le \mathbf{E}^{\text{s}}$ for $\eta_{\text{c}},\eta_{\text{d}}\in (0,1]$.

\begin{lemma}\label{lemma1}
If inputs $\mathbf{P}_{\text{b}}=\mathbf{P}_{\text{c}}-\mathbf{P}_{\text{d}}=\mathbf{P}^{r}_{\text{c}}-\mathbf{P}^{r}_{\text{d}}$ satisfy $\mathbf{P}_{\text{c}}\cdot \mathbf{P}_{\text{d}} = \mathbf{0}$ and $\mathbf{P}_{\text{c}}^\text{r}\cdot \mathbf{P}_{\text{d}}^\text{r}\ge \mathbf{0}$, then $\mathbf{E}^{r}(\mathbf{P}_{\text{c}}^\text{r},\mathbf{P}_{\text{d}}^\text{r})\le \mathbf{E}(\mathbf{P}_{\text{c}},\mathbf{P}_{\text{d}})\le \mathbf{E}^{\text{s}}(\mathbf{P}_\text{b})$.
\end{lemma}

\begin{proof}
First we shall prove that $\mathbf{E}^{r}\le \mathbf{E}$. Subtracting~\eqref{eq:P2} from~\eqref{eq:P1} and substituting the values of $\mathbf{P}_{\text{c}}$ and $\mathbf{P}_{\text{d}}$ from~\eqref{eq:P_subs} and applying basic algebraic operations, we get:
\begin{align}
    & \Delta \mathbf{E}^{r}:=\mathbf{E}-\mathbf{E}^{r}= \mathbf{A}\left(\frac{1}{\eta_{\text{d}}}-\eta_{\text{c}}\right)\min \{\mathbf{P}^{r}_{\text{c}},\mathbf{P}^{r}_{\text{d}}\}\ge 0
\end{align}
This proves the first part of the lemma. % and shows that $\mathbf{E}^{r}\le \mathbf{E}$. 
To prove $\mathbf{E}\le \mathbf{E}^{\text{s}}$, we subtract~\eqref{eq:P1} from~\eqref{eq:P3} and substitute $\mathbf{P}_{\text{b}}=\mathbf{P}_{\text{c}}-\mathbf{P}_{\text{d}}$, which gives
% \begin{align}\label{eq:Delta_B}
%     \Delta \mathbf{{E}^s}=\mathbf{{E}^s}-\mathbf{E}=\Delta t \mathbf{A}(\eta \mathbf{P_{\text{b}}}-\eta_{\text{c}}\mathbf{P_{\text{c}}}+\frac{1}{\eta_{\text{d}}}\mathbf{P_{\text{d}}})
% \end{align}
% Substituting $\mathbf{P_{\text{b}}=P_{\text{c}}-P_{\text{d}}}$ in~\eqref{eq:Delta_B} gives:
% \begin{align}
%     \Delta \mathbf{{E}^s}=&\Delta t \mathbf{A}(\eta \mathbf{P_{\text{c}}}-\eta P_{\text{d}}-\eta_{\text{c}}\mathbf{P_{\text{c}}}+\frac{1}{\eta_{\text{d}}}\mathbf{P_{\text{d}}})\\
\begin{align}
    \Delta \mathbf{E}^{\text{s}}:=\mathbf{E}^{\text{s}}-\mathbf{E}= \mathbf{A}\left[ \mathbf{P}_{\text{c}}(\eta-\eta_{\text{c}})+\mathbf{P}_{\text{d}}\left(\frac{1}{\eta_{\text{d}}}-\eta\right)\right] \ge 0 \label{eq:Delta_B1}
\end{align}
due to $\eta_{\text{c}}\le \eta\le\frac{1}{\eta_{\text{d}}}$.
Thus, 
%proves that $\mathbf{E}^{s}\ge \mathbf{E}$ and shows that 
 $\mathbf{E}^{\text{r}}\le \mathbf{E}\le \mathbf{E}^{\text{s}}$.
\end{proof}

%With $\mathbf{E}^{r}\le \mathbf{E}\le \mathbf{{E}^s}$, 
Lemma~\ref{lemma1} shows that the simplified model overestimates the actual SoC, while the relaxed model underestimates the SoC. Furthermore, from the proof of Lemma~\ref{lemma1}, we can analyze the worst-case SoC model mismatches $\Delta \mathbf{E}^r,\Delta \mathbf{E}^s$. The bounds on the mismatches represent the conservativeness of the two battery models. For the relaxed model, the worst-case mismatch over the trajectory is given by:
  \begin{align}\label{eq:error_relax}
  \Delta \mathbf{E}^{\text{r}}=\mathbf{A}(\frac{1}{\eta_{\text{d}}}-\eta_{\text{c}}) \min\{\mathbf{P}^{\text{r}}_{\text{c}},\mathbf{P}^{\text{r}}_{\text{d}}\}
  \le (\frac{1}{\eta_{\text{d}}}-\eta_{\text{c}}) \mathbf{A} \textbf{1}_T P_\text{max}
 \end{align}
%   &{\color{red} \le  \Delta t \mathbf{A}(\frac{1}{\eta_{\text{d}}}-\eta_{\text{c}}) \max\{\mathbf{P}^{\text{r}}_{\text{c}},\mathbf{P}^{\text{r}}_{\text{d}}\}}\\
%   &\le\Delta t \mathbf{A}(\frac{1}{\eta_{\text{d}}}-\eta_{\text{c}}) \max\{\mathbf{P}_{\text{c}},\mathbf{P}_{\text{d}}\}\\
%   &\le\Delta t \mathbf{A}(\frac{1}{\eta_{\text{d}}}-\eta_{\text{c}}) \max\{\mathbf{0},\mathbf{P}_{\text{b}}\}-\min\{\mathbf{0},\mathbf{P}_{\text{b}}\}.
This worst-case error can further be reduced by including the cutting-plane from~\cite{almassalkhi2014model}: $\mathbf{P}^{r}_{\text{c}}+\mathbf{P}^{r}_{\text{d}} \le \mathbf{1}_T P_{\text{max}}$ in~\eqref{eq:relax_model}, which gives:
\begin{align}\label{eq:B_loss_b}
     \Delta \mathbf{E}^{\text{r}} \le (\frac{1}{\eta_{\text{d}}}-\eta_{\text{c}}) \mathbf{A} \textbf{1}_T \frac{P_{\text{max}}}{2}.
\end{align}
Similarly, the simplified model's mismatch can be written:
  \begin{align}\label{eq:error_simp}
  \Delta \mathbf{E}^{\text{s}}& = \mathbf{A}\left[ (\eta-\eta_{\text{c}})\max\{\mathbf{0},\mathbf{P}_{\text{b}}\}+\left(\eta-\frac{1}{\eta_{\text{d}}}\right)\min\{\mathbf{0},\mathbf{P}_{\text{b}}\}\right].
  \end{align}
% {\color{red} What is key reason for / takeaway from looking at ``worst-case'' bounds? Don't make your paper look like your notes. Motivate and clarify... or remove! }
Note that $\Delta \mathbf{E}^{\text{s}}$ depends on choice of $\eta$.
%we would like to minimize it. As the value of $\mathbf{P_{\text{b}}}$ is not fixed \text{a-priori}, 
Consider the choice of $\eta$ such that $\eta-\eta_{\text{c}}=\frac{1}{\eta_{\text{d}}}-\eta = \frac{1}{2}(\frac{1}{\eta_\text{d}} - \eta_\text{c}) =:\alpha$. Based on this $\eta$, the worst-case simplified model mismatch is
\begin{align}\label{eq:B_loss_c}
    \Delta \mathbf{E}^{\text{s}}= \alpha \mathbf{A} \left[\max\{\mathbf{0},\mathbf{P}_{\text{b}}\}-\min\{\mathbf{0},\mathbf{P}_{\text{b}}\}\right]= \alpha\mathbf{A} |\mathbf{P}_{\text{b}}|\\
\implies \Delta \mathbf{E}^{\text{s}}\le   \alpha\mathbf{A} \mathbf{1}_T P_{\text{max}}= (\frac{1}{\eta_d}-\eta_c)\mathbf{A} \mathbf{1}_T \frac{P_{\text{max}}}{2}.\label{eq:B_loss_d}
\end{align}
From~\eqref{eq:B_loss_b} and~\eqref{eq:B_loss_d}, it can be seen that these worst-case model mismatch bound are equivalent for the given choice of $\eta$. Clearly, for $\eta_\text{c} = 1 = \eta_\text{d}$, both battery models are exact (as is known). However, in practice, $(\frac{1}{\eta_d}-\eta_c)\le 0.2$ for most lithium-based and lead-acid battery technologies (with round-trip efficiencies $>80\%$), which yields model mismatches (well) below $\frac{P_\text{max}}{10} \mathbf{A1}_T$.

Since errors are reasonable, we can employ the relaxed and simplified models as lower and upper bounds, respectively,  in a linear battery optimization formulation that ensures the actual SoC is persistently within SoC limits.

\section{Optimal battery dispatch formulation}\label{sec:prob_form}
Based on the analysis in Lemma~\ref{lemma1}, the two battery models bound the actual SoC. The linear robust battery dispatch (RBD) problem  can then be formulated as follows:
\begin{subequations}\label{eq:opt1}
\begin{align}
   \textbf{(RBD)} \qquad \min_{\mathbf{P}_{\text{c}}-\mathbf{P}_{\text{d}}} & \quad  f(\mathbf{P}_{\text{c}} -\mathbf{P}_{\text{d}})\label{eq:opt1_a}\\
     \textrm{s.t}\quad
    \mathbf{0}\le & \mathbf{1}_TE_0+ \eta_{\text{c}} \mathbf{A}\mathbf{P}_{\text{c}}-\frac{1}{\eta_{\text{d}}} \mathbf{A}\mathbf{P}_{\text{d}} %\qquad | \quad \underline{\lambda}
    \label{eq:opt1_b}\\
    \mathbf{E}_{\text{max}}\ge & \mathbf{1}_T E_0+ \eta \mathbf{A}(\mathbf{P}_{\text{c}}-\mathbf{P}_{\text{d}}) %\qquad | \quad \overline{\lambda}
    \label{eq:opt1_c}\\
    0\le & \mathbf{P}_c\le \mathbf{1}_T P_{\text{max}} %\qquad | \quad \underline{\beta},\overline{\beta}
    \label{eq:opt1_d}\\
     0\le & \mathbf{P}_d\le  \mathbf{1}_T P_{\text{max}} %\qquad | \quad \underline{\gamma},\overline{\gamma}
     \label{eq:opt1_e}\\
     & \mathbf{P}_{\text{c}}+\mathbf{P}_{\text{d}} \le  \mathbf{1}_T P_{\text{max}} %\qquad | \quad \mu
     \label{eq:opt1_f}
\end{align}
\end{subequations}
%where $ \underline{\lambda},\overline{\lambda},\underline{\beta}, \overline{\beta},\underline{\gamma},\overline{\gamma},\mu \in \mathbb{R}^{1\times T}$ are the lagrange multipliers associated with their respective constraints.
\begin{remark}
We can easily adapt~\eqref{eq:opt1} to power systems with $N$ batteries and  modify the objective to $f(\sum_{i=1}^N(\mathbf{P}_{\text{c},i}-\mathbf{P}_{\text{d},i}))$. The formulation can also be augmented by coupling the batteries inputs via power flow equations, e.g.,~\cite{nazir2020optimal}.
\end{remark}
The robust optimization problem in~\eqref{eq:opt1} leads to a conservative battery dispatch. However, the conservativeness is with respect to the objective function. That is, by guaranteeing that the actual SoC trajectory is within its energy limits, the optimization problem always ensures that an optimized power dispatch is realizable. In fact, the conservativeness in the objective depends on the time step width $\Delta t$ and the horizon $T$ (i.e., $\mathbf{A}$) and battery specs ($\eta_\text{c}, \eta_\text{d}, P_\text{max}$).

\begin{remark}
\textcolor{black}{Note that since the results hold for any objective function in~\eqref{eq:opt1_a}, the linear RBD formulation is well-suited in model predictive control (MPC) settings and in unit commitment, security-constrained, and multi-period economic dispatch applications.}
\end{remark}

Next, we illustrate the effectiveness of the proposed approach in~\eqref{eq:opt1} with simulation results.
\section{Simulation Results}\label{sec:simulations}
\begin{figure}[t]
    %\centering
    %\vspace{-9pt}
  \subfloat [\label{fig:batt_track}]{   \includegraphics[width=0.48\linewidth,trim={0 0 0 0.5cm},clip]{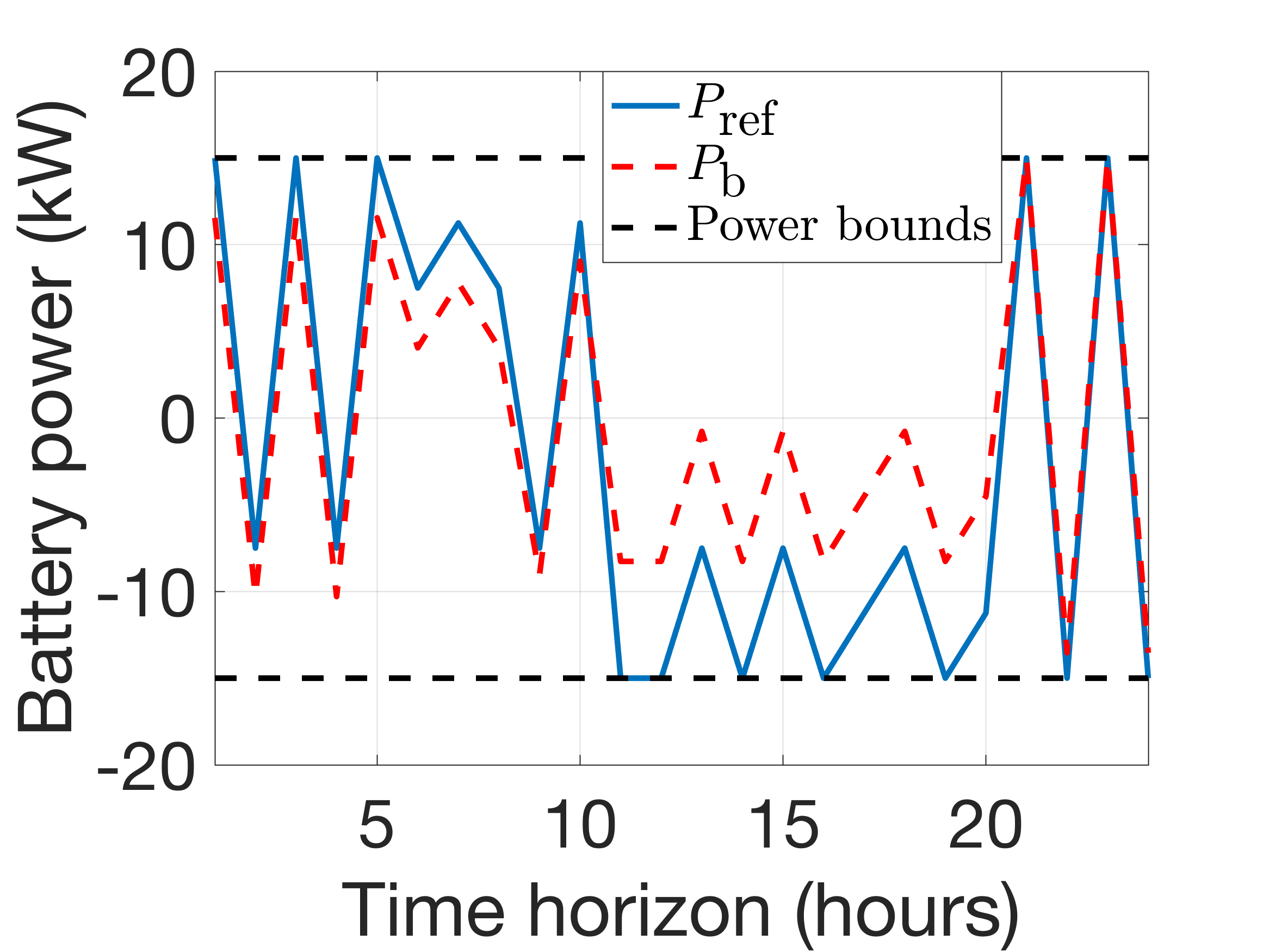}.}
    \hfill
  \subfloat [\label{fig:SoC_limit}]{    \includegraphics[width=0.48\linewidth,trim={0 0 0 0.5cm},clip]{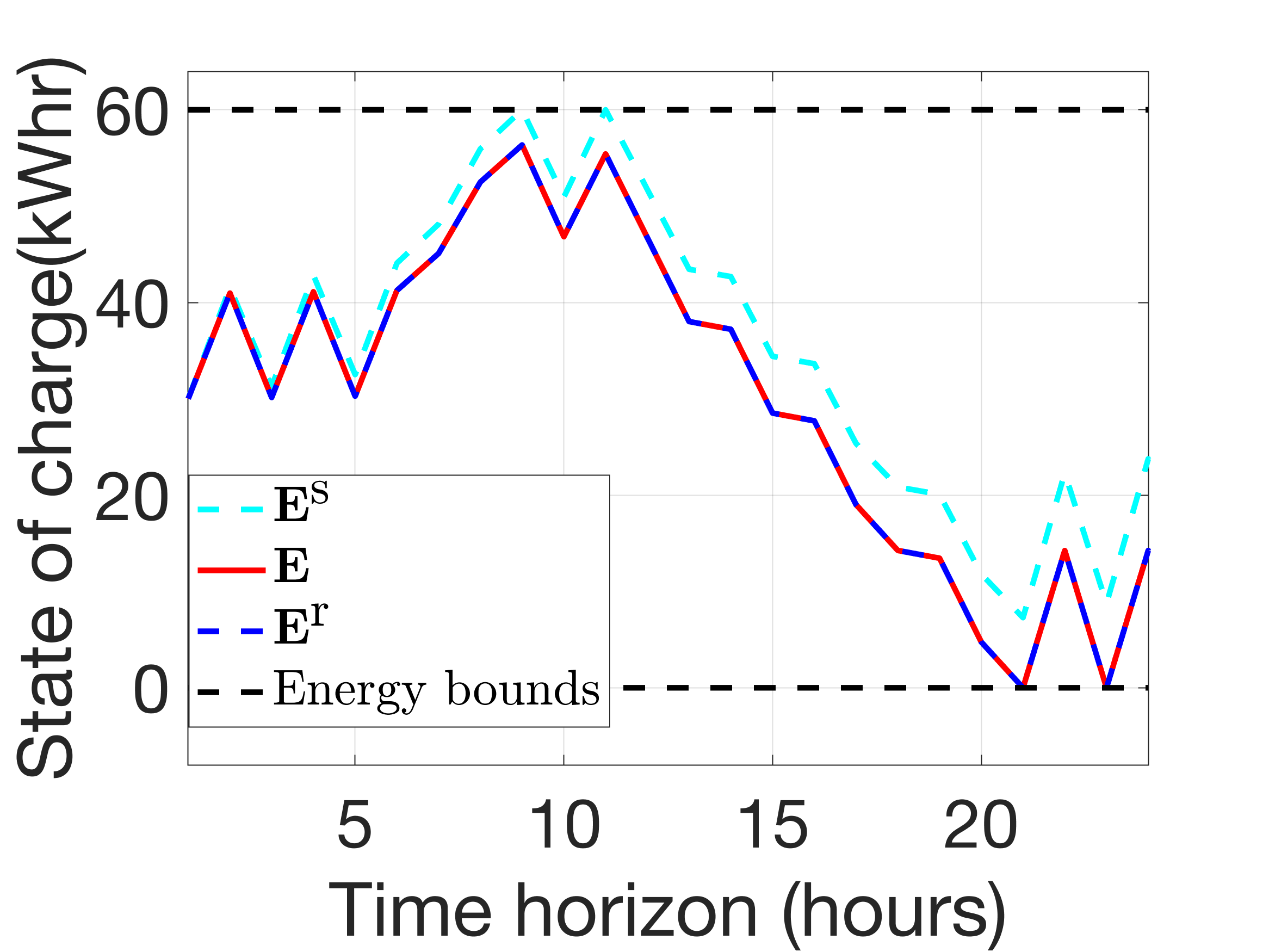}}
  \caption{(a) Tracking a battery reference power signal $P_{\text{ref}}$ with the net battery output $P_{\text{b}}\in [-P_\text{max},P_\text{max}]$.  (b) Comparison between predicted SoC ($\mathbf{E}^{s}, \mathbf{E}^{r}$) and actual SoC $\mathbf{E}$ resulting from optimized dispatch with the energy limits $[0,60]$. Clearly, the actual SoC trajectory $\mathbf{E}$ satisfies energy limits.}
\label{fig:case_study1}
\end{figure}

%To illustrate the effectiveness of the formulation in~\eqref{eq:opt1}, we present simulation results on the optimization problem in. 
Consider a battery with $P_\text{max}=15$kW and $E_\text{max}=60$kWh. Let $\eta_c = 0.95 = \eta_d$ and choose $\eta=(\eta_n + \eta_d)/2 = 1.0013$, which results in  $\alpha=0.0513$. The time-step $\Delta t$ is 1~hour and the control and prediction horizon length $T$ is 24~hours. The objective in~\eqref{eq:opt1_a} is chosen as $\sum_k(P_{\text{ref}}[k]-(P_{\text{ c}}[k]-P_{\text{d}}[k]))^2$.
In Fig.~\ref{fig:batt_track}, the results shows one battery tracking a reference power signal while Fig.~\ref{fig:SoC_limit} compares the predicted (upper and lower bounds) SoC resulting from~\eqref{eq:opt1} to the actual battery SoC obtained from~\eqref{eq:P1}. Fig.~\ref{fig:SoC_limit} illustrates that trajectory $\mathbf{E}$ is within its energy limits, which means the optimized power dispatch $\mathbf{P}_\text{b}$ is guaranteed to be realizable.

Furthermore, to highlight computational efficiency, Table~\ref{table_time} compares the RBD in~\eqref{eq:opt1} to exact mixed-integer (MIP) and non-linear (NLP) formulations as the number of batteries $N$ increases and we track $N P_\text{ref}$. The RBD and MIP are solved using Gurobi~9.1, while the NLP uses IPOPT on a standard laptop. The table shows that the RBD method is 10-200 times faster than MIP for $N\le 200$ batteries. For $N\ge 500$, MIP does not find a solution with MIP-gap $<10\%$ within 3600s. The RBD approach is also 5-50 times faster than the NLP, which only achieves local optimum. Note also that the RBD outperforms the NLP with respect to open-loop tracking performance (RMSE) and is still within 10\% of the globally optimal, exact MIP. The RBD's fast solve time enables a receding-horizon implementation that should greatly reduce RMSE. Thus, with~\eqref{eq:opt1}, we sidestep the challenges with non-convex or integer-based complementarity constraints and provide a linear formulation that guarantees a realizable dispatch.  %However, the conservativeness in the objective is counteracted by the computational benefits.
%Of course, the presented multi-battery test is rather simple and within an OPF setting, we would expect stronger coupling between the battery variables, which should yield greater computational benefits.

\begin{table}[H]
\centering
\caption{Solve time (sec) and power tracking RMSE (kW) comparison with increasing batteries for RBD vs MIP vs NLP}
\begin{tabular}{p{.75cm}p{.75cm}p{1cm}p{.75cm}p{1cm}p{.75cm}p{1cm}}
\toprule
& \multicolumn{2}{p{1.75cm}}{RBD} & \multicolumn{2}{p{1.75cm}}{MIP} & \multicolumn{2}{p{1.75cm}}{NLP}\\
\midrule
{Batteries} & {Time} & {RMSE} & {Time} & {RMSE} & {Time} & {RMSE} \\
\midrule
          %  1 & 1.6 & 1.6 & 0.7\\
            $10$ & $1.7$ & $47.8$ & $16.3$ & $43.7$ & $5.1$ & $54$\\
           $100$ & $3.1$ & $478.7$ & $271.8$ & $437.8$ & $50.5$ & $478.7$\\
           $200$ & $6.3$ & $957.4$ & $1114$ & $866$ & $133.2$ & $1190.2$\\
           $500$ & $11.5$ & $2327.4$ & $-$ & $-$ & $351.6$ & $2415.2$\\
           $1000$ & $22.6$ & $4787.1$ & $-$ & $-$ & $1115$ & $4787.1$\\
\bottomrule
\end{tabular}
\label{table_time}
\end{table}

\subsection{Impact of efficiency}
\begin{figure}[t]
    %\centering
    %\vspace{-9pt}
  \subfloat [\label{fig:SoC_error}]{   \includegraphics[width=0.48\linewidth,trim={0 0 0 0.5cm},clip]{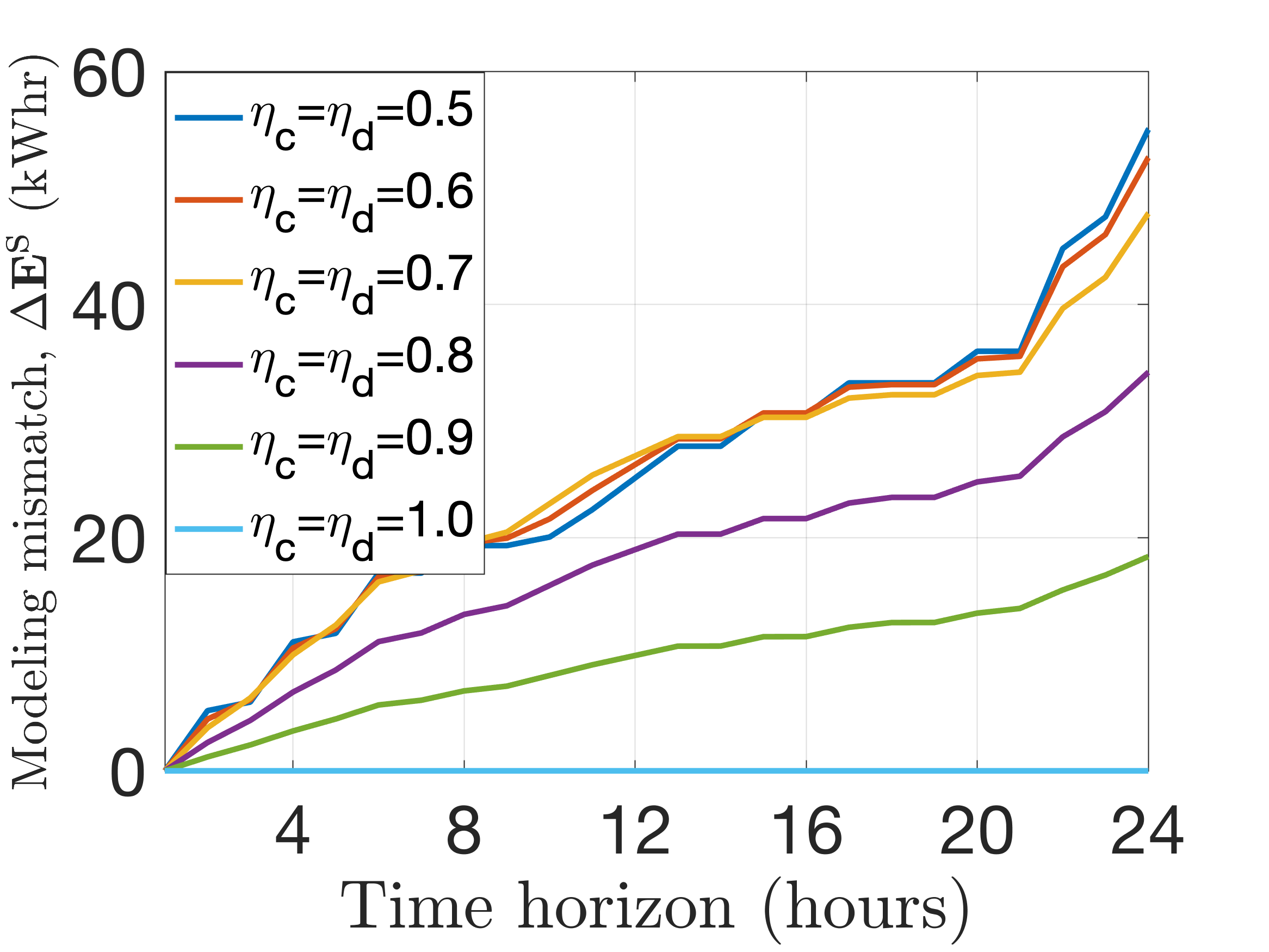}.}
    \hfill
  \subfloat [\label{fig:obj_track}]{    \includegraphics[width=0.48\linewidth,trim={0 0 0 0.5cm},clip]{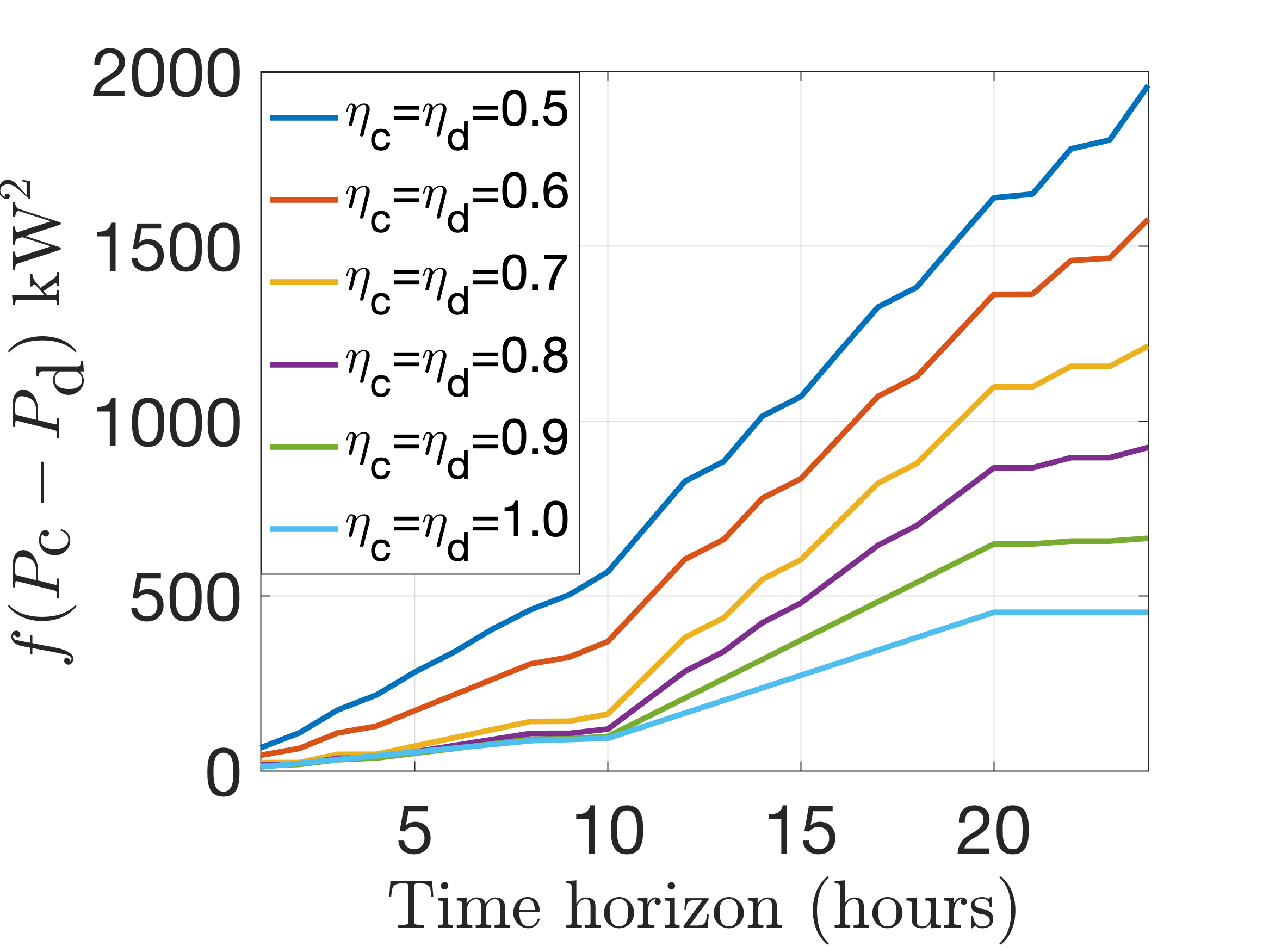}}
  \caption{(a) Modeling mismatch, $\Delta \mathbf{E}^{\text{s}}$ obtained from~\eqref{eq:B_loss_c} for different $\eta_\text{c}= \eta_\text{d}$ efficiencies.  (b) Corresponding cumulative objective function values ($(P_{\text{ref}}[k]-P_{\text{b}}[k])^2$) showing reduced tracking performance with increased modeling mismatch (i.e., lower efficiencies).}
\label{fig:case_study2}
\end{figure}

\textcolor{black}{The model mismatch in the RBD formulation (i.e., $\Delta \mathbf{E^S}$ and $\Delta \mathbf{E^r}$) depends largely on the charge and discharge efficiencies. If the round-trip efficiency is low, then the model mismatch increases, which makes the RBD formulation more conservative. To illustrate the effect of round-trip efficiency on the model mismatch and conservativeness, we repeat the simulations from Fig.~\ref{fig:case_study1}, but over a range of efficiencies. The resulting model mismatch with $\Delta E^S$ is shown in Fig.~\ref{fig:SoC_error} and represents an over-estimate of SoC. The figure illustrates the increased model mismatch as the round-trip efficiency reduces. The corresponding cumulative objective function values are shown in Fig.~\ref{fig:obj_track}, highlighting the reduced tracking performance with lower efficiencies. For applications, such as pumped hydro or hydrogen storage (i.e., electrolyzer+fuel cells), where the round-trip efficiency is around 60\%, the proposed RBD formulation will be conservative and may not be suitable, but is, nonetheless, guaranteed to be \textit{physically realizable}. %That said, pumped hydro and hydrogen storage systems can actually perform simultaneous charging and discharging (unlike conventional batteries), so the inefficiencies due to simultaneous charging and discharging can be quantified in the objective function through energy loss and economic parameters. This means that these specific systems may avoid complementarity constraints by using just the relaxations, which is similar to the formulation presented in~\cite{nazir2020optimal}.
}

% \begin{figure}[t] \label{fig:SoC_error1}
% \centering
% \includegraphics[width=0.4\textwidth]{SoC_error_2.png}
% \caption{Normalized state of charge (SoC) modeling mismatch, $\Delta \mathbf{E}^{\text{s}}$ obtained from~\eqref{eq:B_loss_c} for different $\eta_\text{c}=\eta= \eta_\text{d}$ efficiencies.}
% \end{figure}

\section{Conclusions and Future work}
\label{sec:conclusion}
This paper presented a new linear formulation to optimally dispatch batteries while guaranteeing satisfaction of SoC constraints, without having to resort to a non-convex and/or mixed-integer battery formulations. Through mathematical analysis, we prove that two linear formulations provide upper and lower bounds on the actual SoC, which enables their use as proxy variables in the optimization formulation. Furthermore, we provide worst-case bounds on the conservativeness of this approach. 
These results have the potential to greatly reduce the complexity of energy-constrained battery optimization problems, while guaranteeing satisfaction of actual SoC constraints. 
%This can pave the way for effectively optimizing batteries in power systems with renewable energy variability and electric vehicle charging systems.
Future work will investigate the RBD formulation from~\eqref{eq:opt1} in various MPC and optimal power flow formulations and study the impact of conservativeness in practical applications.

%Bibliography 
\bibliographystyle{IEEEtran}
\small\bibliography{sample.bib}

\begin{IEEEbiography} [{\includegraphics[width=1in,height=1.25in,clip,keepaspectratio]{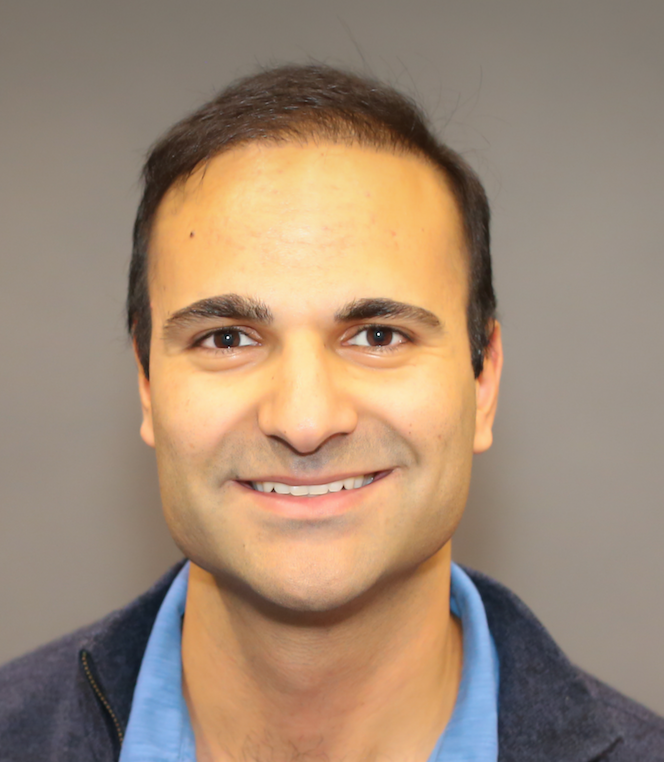}}] {Nawaf Nazir (S’17-M'21)} received the
M.S. degree in electrical engineering from Virginia Polytechnic Institute and State University, Blacksburg, VA, USA, in 2015 and the PhD degree in Electrical Engineering from the Department of Electrical and  Biomedical Engineering at the University of Vermont, Burlington, VT, USA. He is currently a Postdoctoral Researcher at the Pacific Northwest National Lab, Richland, WA, USA.
His research interests include optimization, control and machine learning applied to complex networked systems, and, in particular, emphasizes reliability, resilience and real-time control of energy systems.
\end{IEEEbiography}
\begin{IEEEbiography}[{\includegraphics[width=1in,height=1.25in,clip,keepaspectratio]{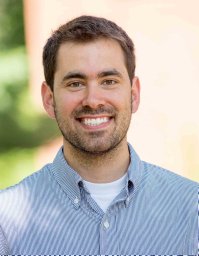}}]{Mads Almassalkhi}(S'06-M'14-SM'19)
	received his B.S. degree in electrical engineering with a dual major in applied mathematics from the University of Cincinnati, Cincinnati, OH, USA, in 2008, and the M.S. degree in Electrical Engineering: Systems from the University of Michigan, Ann Arbor, MI, USA, in 2010, where he also earned his Ph.D. degree in 2013. He is currently an Associate Professor in the Department of Electrical and Biomedical Engineering at the University of Vermont, Burlington, VT, USA. His research interests span multi-timescale control of DERs, energy optimization in power systems, intelligent electrification, and multi-energy systems.
\end{IEEEbiography}

\end{document}